\pdfoutput=1
\documentclass[letterpaper, 10 pt, conference]{ieeeconf}  

\IEEEoverridecommandlockouts                              

\usepackage{graphicx}      
\usepackage{cite}

\usepackage{amssymb}
\usepackage{amsmath}
\usepackage{relsize}
\usepackage{nicefrac,xfrac}
\usepackage{gensymb}
\usepackage{graphicx}
\usepackage{multirow}
\usepackage{subcaption}
\usepackage{algorithm}
\usepackage{algpseudocode}
\usepackage{adjustbox}
\usepackage{nomencl}
\let\labelindent\relax
\usepackage{enumitem}
\usepackage{titlesec}

\title{\LARGE \bf
Explicit Distributed MPC: Reducing Computation and Communication Load by Exploiting Facet Properties*
}

\author{Parth R. Brahmbhatt$^{1}$, Hari S. Ganesh$^{2}$, and Styliani Avraamidou$^{1}$ 
\thanks{*This work is based upon work supported by the National Science Foundation under grant no. CMMI-2328160. Financial support from the  University of Wisconsin-Madison, and the Wisconsin Alumni Research Foundation is also gratefully acknowledged.}
\thanks{$^{1}$Department of Chemical and Biological Engineering, University of Wisconsin-Madison, Madison, WI 53706, USA
        {\tt\small avraamidou@wisc.edu}}%
\thanks{$^{2}$Department of Chemical Engineering, Indian Institute of Technology Gandhinagar, Palaj, Gandhinagar, 382355, Gujarat, India
        {\tt\small hariganesh@iitgn.ac.in}}%
}

\begin{document}
\maketitle
\thispagestyle{empty}
\pagestyle{empty}

\begin{abstract}                
Classical Distributed Model Predictive Control (DiMPC) requires multiple iterations to achieve convergence, leading to high computational and communication burdens. This work focuses on the improvement of an iteration-free distributed MPC methodology that minimizes computational effort and communication load. The aforementioned methodology leverages multiparametric programming to compute explicit control laws offline for each subsystem, enabling real-time control without iterative data exchanges between subsystems. Extending our previous work on iteration-free DiMPC, here we introduce a \textbf{FA}cet-based \textbf{C}ritical region \textbf{E}xploration \textbf{T}echnique for iteration-free DiMPC (FACET-DiMPC) that further reduces computational complexity by leveraging facet properties to do targeted critical region exploration. Simulation results demonstrate that the developed method achieves comparable control performance to centralized methods, while significantly reducing communication overhead and computation time. In particular, the proposed methodology offers substantial efficiency gains in terms of the average computation time reduction of 98\% compared to classic iterative DiMPC methods and 42\% compared to iteration-free DiMPC methods, making it well-suited for real-time control applications with tight latency and computation constraints.
\end{abstract}


\section{Introduction} 

Model Predictive Control (MPC) is an optimal control technique widely used since the late 1980s in industries like chemical plants and oil refineries \cite{ qin1997overview}. Despite its ability to handle constraints and optimize process inputs, the adoption of MPC for large-scale systems faces challenges due to increased computational burdens\cite{trodden2017distributed}. Centralized MPC (CMPC) calculates control inputs by solving an optimal control problem (OCP) for the entire plant at each time step, but its high computational cost in the case of plants with multiple control variables can hinder real-time performance. In contrast, decentralized MPC (DeMPC) reduces complexity by decoupling the plant into subsystems, decreasing computational effort. However, the decoupling of subsystems limits the ability to capture interactions between them, leading to suboptimal control and potential instabilities \cite{cui2002performance, diangelakis2016decentralized}.

Distributed MPC (DiMPC) offers a balance by combining the decoupled architecture of DeMPC with CMPC's ability to capture interactions between subsystems. The DiMPC architecture can be categorized into the non-cooperative and cooperative forms\cite{negenborn2014distributed, farina2012distributed}. Cooperative DiMPC optimizes control inputs considering both local information and global system performance\cite{stewart2010cooperative}. However, it involves significant communication and computational costs due to multiple OCPs and intermediate iterations. 

Recent advancements focus on improving DiMPC's efficiency. Multiparametric (mp) programming-based methods reduce online computation by generating offline control laws as affine functions \cite{pistikopoulos2020multi}. Despite these improvements, communication burdens remain \cite{saini2023noncooperative}. To address these challenges, iteration-free approaches have been developed ~\cite{camponogara2002distributed,chen2016improved}. One such method, iteration-free cooperative DiMPC (IF-mpDiMPC), reduces communication loads by solving precomputed MPC control policy functions, minimizing delays, and enhancing robustness against data loss and cybersecurity. This can lead to an unnecessarily large search space, as regions sharing a hyperplane may not share a common facet. 

In this work, we introduce a FACET-DiMPC: \textbf{FA}cet-based \textbf{C}ritical region \textbf{E}xploration \textbf{T}echnique for iteration-free \textbf{DiMPC}, that leverages the facet properties of explicit controllers to enhance computational efficiency and reduce communication overhead. Building on our previous work on iteration-free DiMPC in Saini \textit{et al.}\cite{saini2025iteration}, we propose refining the exploration of critical regions. In Saini \textit{et al.}\cite{saini2025iteration}, we explored the current critical region where past inputs were captured, along with other regions that share a hyperplane with the current region (section 2.3 presents more details on this method). In this work, we aim to improve the computational efficiency of our previous work on iteration-free DiMPC~\cite{saini2025iteration} by considering regions that share facets, further refining the search for optimal control inputs. Closed-loop numerical simulations on randomly generated plants demonstrate that the proposed iteration-free DiMPC method outperforms that reported in our previous work~\cite{saini2025iteration}, making it highly effective for real-time distributed control applications where reducing latency and communication are crucial. 

\section{Background on Iterative \& Iteration-Free DiMPC using Multi-Parametric Programming}
\label{sec:background}

This section provides an overview of the formulation and solution methods for multiparametric programming MPC (mpMPC) \cite{pistikopoulos2020multi,pappas2021multiparametric}, Iterative mpDiMPC (I-mpDiMPC) \cite{saini2023noncooperative}, and our previous work on iteration-free mpDiMPC (IF-mpDiMPC) \cite{saini2025iteration}, respectively.

\subsection{Multi-Parametric (mp) MPC (Explicit-MPC)}
MPC is an optimization-based control strategy that derives control actions by solving an Optimal Control Problem (OCP) shown below at each control time step:  
\begin{equation}
\begin{aligned}
   \min_{\textbf{U}(k)} &  J(\textbf{x}(k), \textbf{U}(k))\\
    \text{s.t.} & ~ \textbf{x}(l+1|k) = A\textbf{x}(l|k) + B\textbf{u}(l|k) , \forall ~ l \in \mathbb{I}_{0:N_p-1} \\
    & \textbf{x}(l|k) \in \chi \ \forall \ l \in \mathbb{I}_{1:N_p} \\
    & \textbf{x}(N_p|k) \in \Omega \\
    & \textbf{u}(l|k) \in \mathcal{U} \ \forall \ l \in \mathbb{I}_{0:N_p-1} \\
    & \textbf{x}(0|k) = \textbf{x}(k)
\end{aligned}
\label{eq: CMPCFormulation}
\end{equation}

\noindent where, $J$ is the objective function of OCP. $\textbf{x} \in \mathbb{R}^{n_x}$ is the state vector, $\textbf{u} \in \mathbb{R}^{n_u}$ is the input vector, $A \in \mathbb{R}^{n_x \times n_x}$ is the system matrix, $B \in \mathbb{R}^{n_x \times n_u}$ is the input matrix, and \textit{k} is the sample time instant. $\chi$, $\Omega$, and $\mathcal{U}$ are the bounding sets for state vector, terminal state vector, and control input vector, respectively, $\textbf{U}(k) = [\textbf{u}(0|k), \textbf{u}(1|k), ..., \textbf{u}(N_p-1|k)$ is the control input decision variable vector, $j|k$ is the $j$ steps ahead prediction from the current time $k$, $l$ is a time step related index, and $N_p$ is the prediction horizon. Solving equation~\eqref{eq: CMPCFormulation} results in the optimal decision vector $\textbf{U}^*(k)$, among which, the first control input $\textbf{u}^*(k) \ (= \textbf{u}(0|k))$ is implemented in the system.

The MPC described in equation~\eqref{eq: CMPCFormulation} is also referred to as online MPC, as the OCP is solved at each sample time step $k$ to obtain the optimized control inputs. The same results can be obtained by determining the explicit expressions of the control law offline through the mp programming approach, reducing online calculations to point location and function evaluation \cite{oberdieck2016multi}. In the MPC formulation, the state vector $\textbf{x}(k)$ is usually the only uncertainty (unknown parameter) that becomes known online at each sample time instant. However, variables such as the previous control input information, measured output information, varying set points, etc., could also be considered as unknown parameters if they are made available only online during controller deployment. Defining the set of the unknown parametric variables as $\theta$ vector ($\theta = \textbf{x}(k)$) and substituting in eq.~\eqref{eq: CMPCFormulation} to reformulate the OCP into,
\begin{equation}
\begin{aligned}
    \min_{\textbf{U}} \frac{1}{2}\textbf{U}^TH\textbf{U} + (\theta ^TH_{t} + c)\textbf{U} \\
    G\textbf{U} \leq b + F\theta
\end{aligned}
\label{eq: multiparametricProblem}
\end{equation}

\noindent where, $H \in \mathbb{R}^{n_u \times n_u}$, $H_{t} \in \mathbb{R}^{n_x \times n_u}$, $c \in \mathbb{R}^{n_u}$, $G \in \mathbb{R}^{n_c \times n_u}$, $b \in \mathbb{R}^{n_c}$,$F \in \mathbb{R}^{n_c \times n_x}$ , and $n_{c}$ is the number of constrains.
The problem can be solved multiparametrically to express the optimal control inputs as explicit affine functions of vector $\theta$, which remain both valid and optimal in a specific polyhedral space or critical region. The explicit (multiparametric) solution to eq.~\eqref{eq: multiparametricProblem} can be expressed mathematically as,
\begin{equation}
 \begin{aligned}
     \textbf{U}^* = f^v(\theta) ~~ if ~~ CR^v : \Phi^v\theta \leq \phi^v \\  
     \forall ~~ v = 1, 2, 3, ..., n_{CR}
 \end{aligned}
 \label{eq: mpMPCSol}
 \end{equation}

 \noindent where $f$ represents the affine function defined for each critical region $v$ with $\Phi \ \text{and} \ \phi$ as the corresponding inequality matrices and $n_{CR}$ represents the total number of critical regions. 

\subsection{Iterative Explicit DiMPC (I-mpDiMPC)}

In Distributed MPC (DiMPC), the system is divided into $M$ interacting subsystems, with each subsystem $i$ having its own local controller that considers interaction dynamics. These controllers exchange information and iteratively solve a plant-wide objective to determine control inputs. At each sample time, local controllers calculate their inputs through an iterative process, using the latest information from others. This continues until a termination criterion is met, after which the final control input is implemented. At the next sample time, this process repeats using the previous control inputs as the starting point. The state dynamics of the $i^{th}$ subsystem are:
\begin{equation}
    \textbf{x}_i(k+1) = A_i\textbf{x}_i(k) + \sum_{j=1}^M B_{i,j}\textbf{u}_j(k)
\end{equation}

\noindent where, $\textbf{x}_i \in \mathbb{R}^{n_{x,i}}$, $ A \in \mathbb{R}^{n_{x,i} \times n_{x,i}}$, $B_{i,j} \in \mathbb{R}^{n_{x,i} \times n_{u,j}}$, and $\textbf{u}_j \in \mathbb{R}^{n_{u,j}}$. The overall system can be assembled as,
\begin{equation}
    \textbf{x}(k+1) = A\textbf{x}(k) + \sum_{j=1}^M B_{j}\textbf{u}_j(k)
    \label{eq: overall_system_model}
\end{equation}
\noindent
\noindent
where $\textbf{x}(k) = [\textbf{x}_1(k), \textbf{x}_2(k), \ldots, \textbf{x}_M(k)] \in \mathbb{R}^{n_x}$, where $n_x = \sum_{i=1}^M n_{x,i}$. Here $A$ is a block-diagonal matrix: $A = \text{diag}(A_1, A_2, \ldots, A_M)$. The matrix $B_j$ is given by $B_j = [B_{1,j}^T, B_{2,j}^T, \ldots, B_{M,j}^T]^T$.

In the I-mpDiMPC approach, each local controller in the distributed system solves its OCP using multiparametric programming. I-mpDiMPC retains the iterative communication process between local controllers. The control laws are precomputed offline as in mpMPC, but the solution at each time step still relies on exchanging information between the controllers to ensure convergence.

Each local controller solves the following OCP iteratively:
\begin{equation}
\min_{U_i(k)} J(x(k), U(k)) = \sum_{i=1}^{M} \rho_i J_i(x_i(k), U_i(k))
\end{equation}
subject to the subsystem dynamics and constraints. Here, $J_i$ is the cost function for each local controller $i$, and $\rho_i$ are weighting factors for the global objective function. During optimization, the control inputs of other local controllers, $\textbf{U}_j(k)$, $\ \forall \ j \in \mathbb{I}_{1:M \backslash i}$, remain unchanged. Once optimal control inputs for all $M$ controllers at the $p^{th}$ intermediate iteration are obtained, the next iterate is calculated using a convex combination of the current and previous iterates:
\begin{equation}
    \textbf{U}^{(p)}(k) = \textbf{w}\textbf{U}^{(p-1)}(k) + (1-\textbf{w})\textbf{U}^{(p)}(k)
\end{equation}

The weights $\textbf{w}$ influence convergence. An adaptive algorithm, as described in \cite{wegstein1958accelerating}, can be used to accelerate convergence by adjusting $\textbf{w}$ at each step based on previous results. The vector $\textbf{w} = [ w_1, w_2, ..., w_i, ..., w_m ]$ contains weights for each subsystem $i$, ensuring $w_i + (1 - w_i) = 1$ for each subsystem. The iterative procedure continues until a predefined convergence criterion is satisfied, typically when the difference between the control inputs in consecutive iterations falls below a certain threshold.

An explicit version of DiMPC was developed by \cite{saini2023noncooperative}, expressing future states as a function of the current state $\textbf{x}(k)$ and future control actions:
\begin{equation}
    \textbf{x}(l+1|k) = A^{l+1}\textbf{x}(k) + \mathlarger{\sum}_{q=0}^{l} A^q\sum_{j=1}^MB_j\textbf{u}_j(l-q|k)
\label{eq: DiMPC_futureStatesUpdateEquation}
\end{equation}

For each local controller $i$, the control actions of other controllers $\textbf{U}j$ are treated as unknown parameters, grouped into $\theta_i = [ \textbf{x}(k), \textbf{U}1, ..., \textbf{U}{i-1}, \textbf{U}{i+1}, ..., \textbf{U}_M]$. Substituting this into the optimization problem results in:
\begin{equation}
\begin{aligned}
    \mathcal{F}_i^{(p)}(\theta_i) = \min_{\textbf{U}_i} \frac{1}{2}\textbf{U}_i^TH_i\textbf{U}_i + (\theta_i^TH_{t,i} + c_i)\textbf{U}_i \\
    G_i\textbf{U}_i \leq b_i + F_i\theta_i
\end{aligned}
\label{eq: DiMPC_multiparametricProblem}
\end{equation}

This problem can be solved multiparametrically, yielding optimal control inputs as explicit affine functions of $\theta_i$ within specific critical regions:
\begin{equation}
 \begin{aligned}
     \textbf{U}_i^{(p)} = f_i^v(\theta_i) ~~ if ~~ CR_i^v : \Phi_i^v\theta_i \leq \phi_i^v \\  
     \forall ~~ v = 1, 2, 3, ..., n_{CR,i}
 \end{aligned}
 \label{eq: M1_DiMPCSol}
 \end{equation}
Here, $f_i$ is an affine function for each critical region. 

The I-mpDiMPC approach balances the reduction in online computation achieved by mp programming with the need for communication and iteration to ensure a globally optimal solution. While the use of precomputed control laws reduces the time required to solve each iteration, the communication overhead and the iterative nature of the algorithm can still introduce delays, particularly in real-time control systems with large numbers of subsystems.

\subsection{Iteration-Free Explicit DiMPC (IF-mpDiMPC)}
The I-mpDiMPC algorithm was improved to be iteration-free by simultaneously solving offline-generated control law expressions for local controllers during online deployment, eliminating the need for iterations and critical region searches at each time step \cite{saini2025iteration}. In the local IF-mpDiMPC controller $i$, the parametric vector $\theta_i$ includes the state information $\textbf{x}(k)$ and control inputs from other controllers, represented as:
\begin{equation}
    \theta_i = \big[\underbrace{\textbf{x}(k)}_\text{$\Bar{\theta}$}, \ \underbrace{\textbf{U}_1, \textbf{U}_2, ..., \textbf{U}_{i-1}, \textbf{U}_{i+1}, ... \textbf{U}_M}_\text{$\Bar{\textbf{V}}_i$} \ \big]
\end{equation}

\noindent where, $\Bar{\textbf{V}}_i$ is the control input information of all other systems and $\Bar{\theta}$ is the remaining unknown parameters. In this study, $\Bar{\theta} = \textbf{x}(k)$. Substituting $\theta_i$ in eq.~\eqref{eq: M1_DiMPCSol} yields the expression of the control law given by,

 \begin{equation}
 \begin{aligned}
     \textbf{U}_i = f_i^v(\Bar{\theta}, \Bar{\textbf{V}}_i) ~~ if ~~ CR_i^v : \Phi_{i,1}^v\Bar{\theta} + \Phi_{i,2}^v\Bar{\textbf{V}}_i \leq \phi_i^v \\
     \forall ~~ v = 1, 2, 3, ..., n_{CR,i}
 \end{aligned}
 \label{eq: M1_DiMPCSol2}
 \end{equation}

\noindent where, $\Phi_{i}^v = \begin{bmatrix}
    \Phi_{i,1}^v & \textbf{0} \\ \textbf{0}  & \Phi_{i,2}^v
\end{bmatrix}$ is a block-diagonal matrix. Similarly, multiparametric solutions for all the $M$ local controllers can be determined, and the explicit expressions of the control law can be assembled as given below,
\begin{equation}
\begin{aligned}
    & \textbf{U} = g_1^v(\Bar{\theta}) \quad \text{if} \quad CR_1^v : \Phi_{1,1}^v\Bar{\theta}_1 + \Phi_{1,2}^v\Bar{\textbf{V}}_1 \leq \phi_1^v, \\
    & \qquad \forall ~ v = 1, 2, ..., n_{CR,1} \\[0pt]
    & \vdots \\[0pt]
    & \textbf{U} = g_M^v(\Bar{\theta}) \quad \text{if} \quad CR_M^v : \Phi_{M,1}^v\Bar{\theta}_M + \Phi_{M,2}^v\Bar{\textbf{V}}_M \leq \phi_M^v, \\
    & \qquad \forall ~ v = 1, 2, ..., n_{CR,M}
\end{aligned}
\label{eq: M1_DiMPCSol4}
\end{equation}

\noindent where, $\textbf{U}$ is composed of  $\Bar{\textbf{V}}_i$ and $\textbf{U}_i$. $\Bar{\textbf{V}}_i$ is a subset of $\textbf{U}$ that includes the control inputs for all subsystems except the $i^{th}$ one, i.e., $[\textbf{U}_1, ..., \textbf{U}_{i-1}, \textbf{U}_{i+1}, ..., \textbf{U}_M]$.
This approach, while making the method iteration-free, still requires solving equations for all combinations of critical regions. This process can be computationally intensive, as it requires solving equations $n_{CR,1} \times n_{CR,2} \times \ldots \times n_{CR,M}$ times. To address this complexity, the updated version of IF-mpDiMPC algorithm was developed by reducing the critical region combinations and frequency of solving these equations online at each sample time.

Instead of searching all possible regions at each time step, only regions that share a hyperplane with the current optimal region are considered. Using $\theta_i$, the optimal critical region is identified. At time $k+1$, this region and its nearest neighbors, which share common hyper planes are considered for control input calculation, reducing the set size compared to standard IF-mpDiMPC. This is a relatively trivial task, as regions that share a common hyperplane can be easily identified. To find the neighbors sharing the critical regions, we compare through all the combination of given pair of half-space hyper-planes in critical regions. As shown in Figure \ref{fig:neighbors}, $R1$ and $R2$ critical regions are sharing the common hyperplane $H1$ for all the cases $a$, $b$, and $c$. However, the challenge arises because regions that share a common hyperplane do not always share a common facet, which means more critical region (CR) combinations need to be solved in subsequent computations in IF-mpDiMPC. This increases computational complexity and can result in more regions to handle.

\section{FACET-DiMPC}

\begin{figure}[htbp]
    \centering
    \begin{subfigure}[b]{0.155\textwidth}
        \centering
        \includegraphics[width=\textwidth]{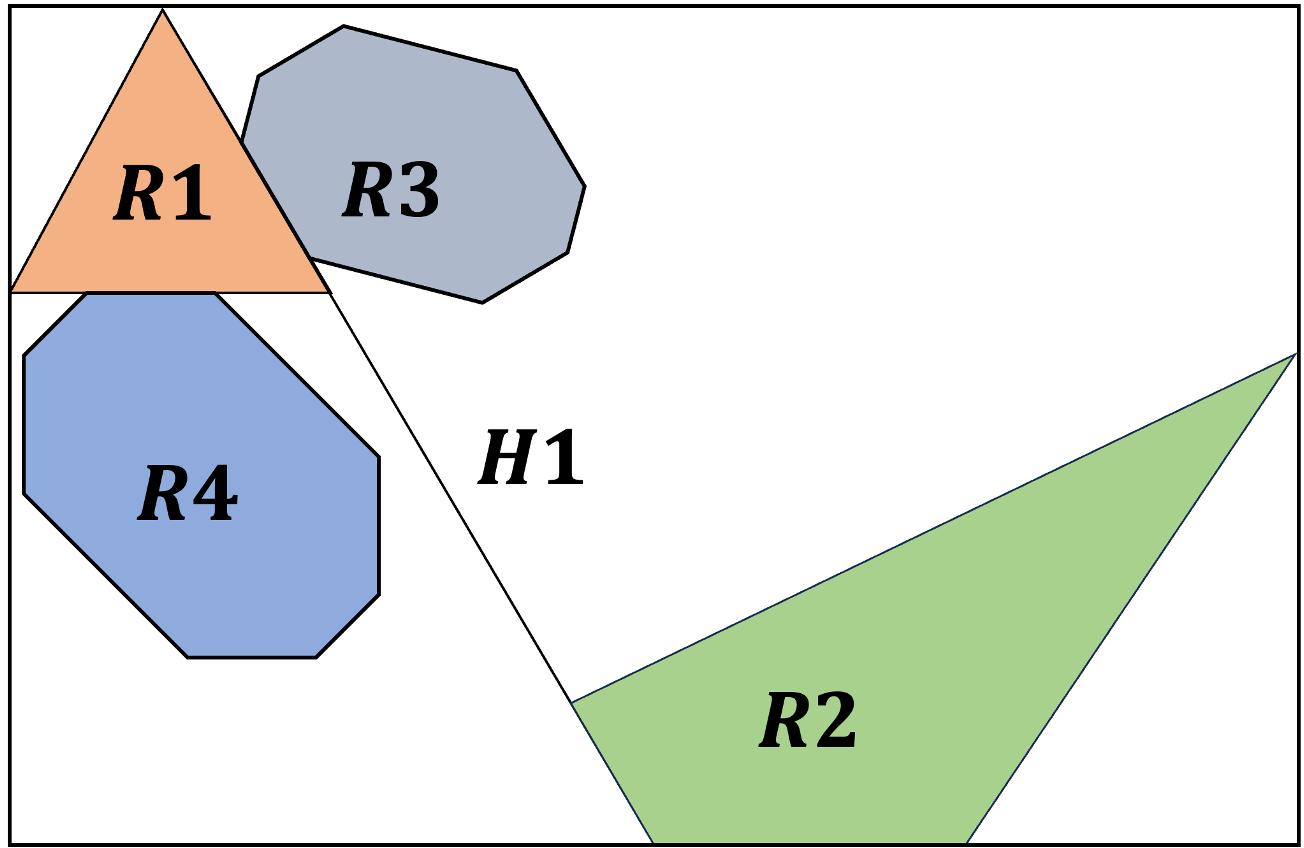}
        \caption{}
        \label{fig:subfig1}
    \end{subfigure}
    \begin{subfigure}[b]{0.155\textwidth}
        \centering
        \includegraphics[width=\textwidth]{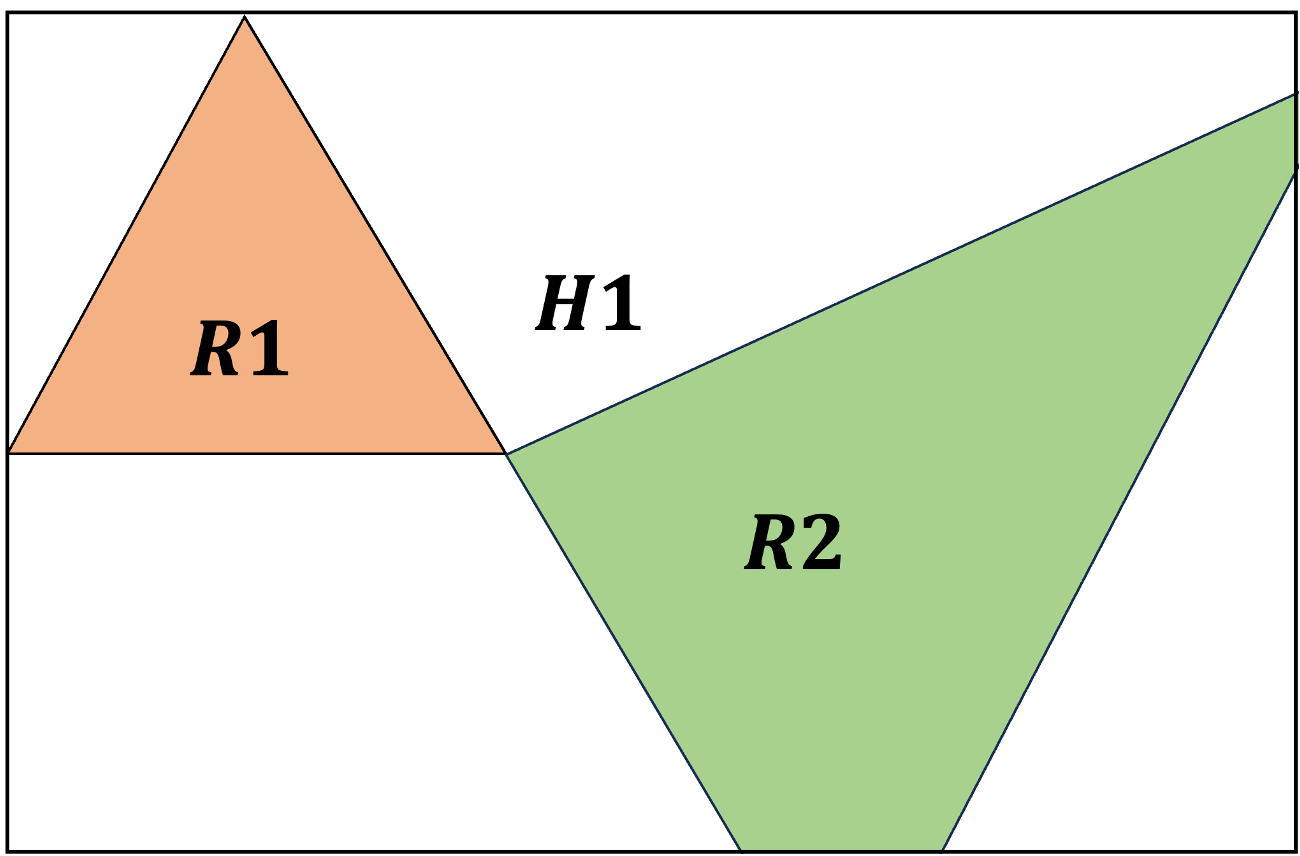}
        \caption{}
        \label{fig:subfig2}
    \end{subfigure}
    \begin{subfigure}[b]{0.155\textwidth}
        \centering
        \includegraphics[width=\textwidth]{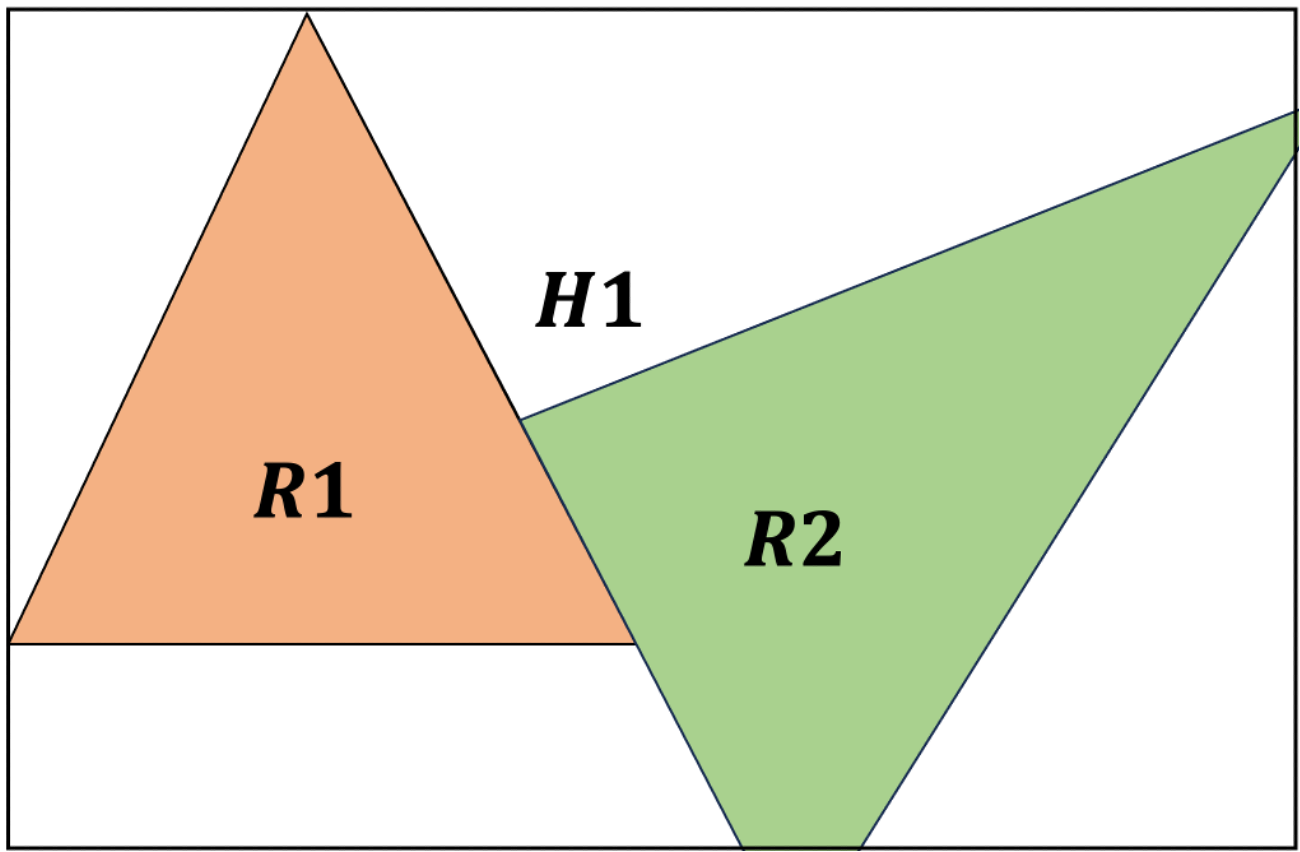}
        \caption{}
        \label{fig:subfig3}
    \end{subfigure}
    
    \caption{Possibilities exist for regions $R1$ and $R2$ sharing a common hyperplane $H1$.}
    \label{fig:neighbors}
\end{figure}

In this section, we present a new facet-based approach aimed at reducing the exploration space of critical regions for the development of the \textit{FACET-DiMPC}. Unlike our previous method, which relied on comparing common hyperplanes to identify neighboring regions, this work advances the process of identifying neighboring regions by solving a linear programming (LP) problem (Eq. \ref{eq: neighborLP}) to rigorously detect region pairs sharing a common facet. Crucially, only in case $c$ of Figure \ref{fig:neighbors}, regions $R1$ and $R2$ share a common facet $H1$. Identifying such true adjacencies is a non-trivial task, unlike simply finding shared hyperplanes. However, this rigorous identification is the key to the FACET-DiMPC framework's efficiency, as it significantly reduces the number of critical region combinations that must be solved, thereby decreasing the online computational effort.

Determining neighboring regions in mp solution space is complex because regions sharing a common hyperplane may not always share a common facet. This complexity arises from various geometric relationships that can exist between regions (Figure \ref{fig:neighbors}):
\begin{enumerate}
    \item \textbf{Non-intersecting Regions}: Two regions may share a common hyperplane but do not intersect at any point (Figure \ref{fig:subfig1}).
    \item \textbf{Point Intersection}: Regions may intersect at a single point on the shared hyperplane (Figure \ref{fig:subfig2}).
    \item \textbf{Common Facet}: Regions genuinely share a facet represented by the common hyperplane (Figure \ref{fig:subfig3}).
\end{enumerate}

To determine if region \(R_i\) shares its \(j^{\text{th}}\) facet, represented by hyperplane \(H_j^i\), with region \(R_k\), we first check if \(H_j^i\) is also a facet of \(R_k\). If it is, we solve an optimization problem to determine whether \(H_j^i\) is a common facet \cite{airan2013linear} :
\begin{equation}
\begin{aligned}
& \text{if } (R_i \cap R_k) \neq \varnothing, \\
& \{\max_{t,x} ~t, \\
\text{s.t.} ~&  x \in R_i, \\
& x \in R_k, \\
& t \leq d(x, H_l^i) \quad \forall l \in \{1, \ldots, n_{f_i}\} - \{j\}\}
\end{aligned}
\label{eq: neighborLP}
\end{equation}
Here, \(t\) is the minimum distance between point \(x\) and all facets of region \(R_i\) except the \(j^{\text{th}}\) facet. The distance function \(d(x, H_l^i)\) calculates the distance from point \(x\) to the \(l^{\text{th}}\) facet of region \(R_i\). The solution to the optimization problem (\ref{eq: neighborLP}) can be:

\begin{itemize}
    \item Infeasible: If regions do not share a facet.
    \item Feasible with \(t = 0\): Regions intersect at a point but do not share a facet.
    \item Feasible with \(t > 0\): Regions share a common facet.
\end{itemize}

After identifying and storing these nearest-neighbor relationships offline, the FACET-DiMPC approach is applied at every time step. If the optimal control input is not found within the pre-computed nearest neighbors at a sample time, the algorithm reverts to the iterative I-mpDiMPC approach. This ensures that the system always has a valid control input, even if additional iterations are required in certain cases, and the worst-case execution time of FACET-DiMPC matches that of I-mpDiMPC. The control law expressions for both I-mpDiMPC and FACET-DiMPC are identical, though FACET-DiMPC reformulates them into simultaneous equations. This reformulation doesn't change the solution space or optimal solutions, so the stability properties of the iterative algorithms apply to FACET-DiMPC as well, eliminating the need for additional stability proofs.


\section{Results and Discussion}

In this section, we present the performance analysis of the proposed FACET-DiMPC and compared with previous approaches shown in section-\ref{sec:background}. The algorithm is evaluated through numerical simulations on a range of systems with varying numbers of subsystems, demonstrating the trade-offs between computational efficiency and communication load.
\subsection{Case Studies}
To evaluate the performance of the proposed iteration-free algorithms, systems consisting of $M \in \{2, 3, 4, 5\}$ subsystems were used. Each subsystem includes a plant with two states and one input $(n_{x,i} = 2, n_{u,i} = 1, \forall i \in 1:M)$. Elements of the state matrix $A_i$ and input matrix $B_{i,j}$ are randomly chosen from $[-1, 1]$, with subsystem coupling reflected in the structure of $B_{i,j}$. State bounds are randomly chosen from $[-100, -10]$ (lower) and $[10, 100]$ (upper), while input bounds are selected from $[-5, -1]$ (lower) and $[1, 5]$ (upper). Both prediction and control horizons are set to $N_p = 3$. Following the methodology in \cite{wang2022improved}, 100 random stable and controllable plants were generated for each $M$ to evaluate computational performance across diverse system dynamics and controller configurations. The same initial conditions are used across all control architectures for consistency.

\begin{figure}[htbp]
    \centering
    \includegraphics[width=0.6\columnwidth]{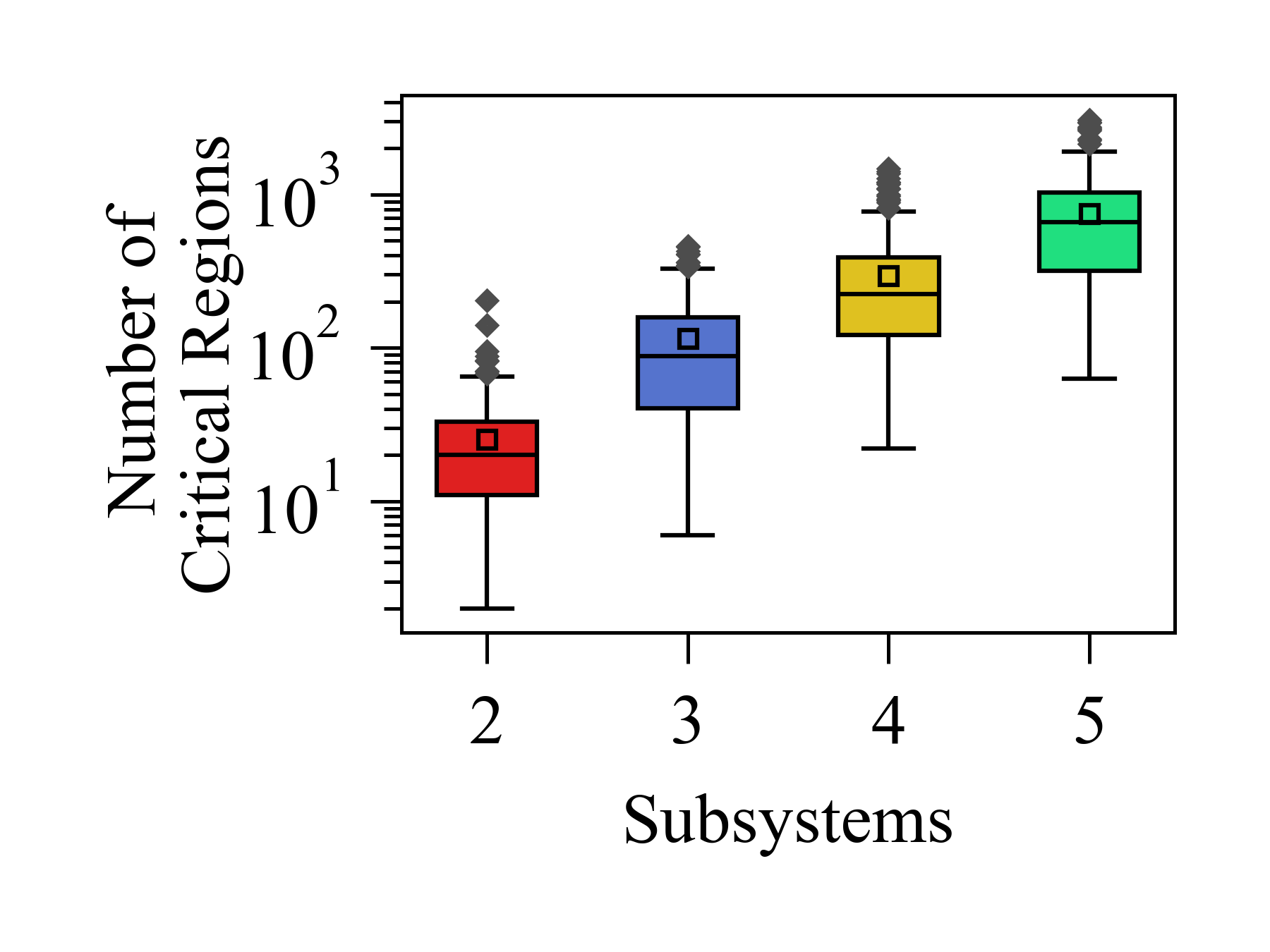}
    \caption{Number of critical regions for different subsystems.}
    \label{fig:critical_regions}
\end{figure}

For instance, for $M = 2$, a sample plant is:

\begin{equation}
    \begin{aligned}
        \textbf{x}_1(k+1) = A_1\textbf{x}_1(k) + B_{1,1}u_1(k) + B_{1,2}u_2(k) \\
        \textbf{x}_2(k+1) = A_2\textbf{x}_2(k) + B_{2,1}u_1(k) + B_{2,2}u_2(k) 
    \end{aligned}
\end{equation}

\noindent where, $\textbf{x}_1 \in \mathbb{R}^2$, $\textbf{x}_2 \in \mathbb{R}^2$,
    $A_1 = \begin{bmatrix}
    0.1645  &   0.7399  \\
    0.0815  &   -0.4704 
    \end{bmatrix}$, $B_{1,1} = \begin{bmatrix}
    -0.3639 \\  -0.7616
    \end{bmatrix}$, $B_{1,2} = \begin{bmatrix}
    0.8797  \\  0.2911 
    \end{bmatrix}$, $A_2 = \begin{bmatrix}
    -0.0411  &   0.0894  \\
    0.2786  &   0.2946 
    \end{bmatrix}$, $B_{2,1} = \begin{bmatrix}
    0.0878  \\  0.4421 
    \end{bmatrix}$, and $B_{2,2} = \begin{bmatrix}
    0.0450  \\ 0.9874 
    \end{bmatrix}$. State and input bounds for this example are:
    \begin{equation}
        \begin{aligned}
            \begin{bmatrix}
                -63.5878    \\   -59.6464
            \end{bmatrix} \leq  \textbf{x}_1 \leq \begin{bmatrix}
                    29.6809 \\ 19.5218
            \end{bmatrix}
            -1.2686 \leq u_1 \leq 3.5116\\
            \begin{bmatrix}
                -67.0765    \\  -31.2846
            \end{bmatrix} \leq  \textbf{x}_2 \leq \begin{bmatrix}
                19.8728 \\  15.7232
            \end{bmatrix}
        -1.1090 \leq u_2 \leq 4.0879\\
        \end{aligned}
    \end{equation}
Explicit solutions and critical region neighbors for each local controller were precomputed and stored for online use. Fig.~\ref{fig:critical_regions} shows the total number of critical regions across all local controllers. The number of critical regions, which increases exponentially with the number of subsystems, remains constant across all three mp programming approaches, with differences only in how optimal control inputs are computed.


\subsection{Control Performance}
\begin{figure}
  \centering
        \begin{subfigure}[b]{0.23\textwidth}
            \includegraphics[width=\linewidth]{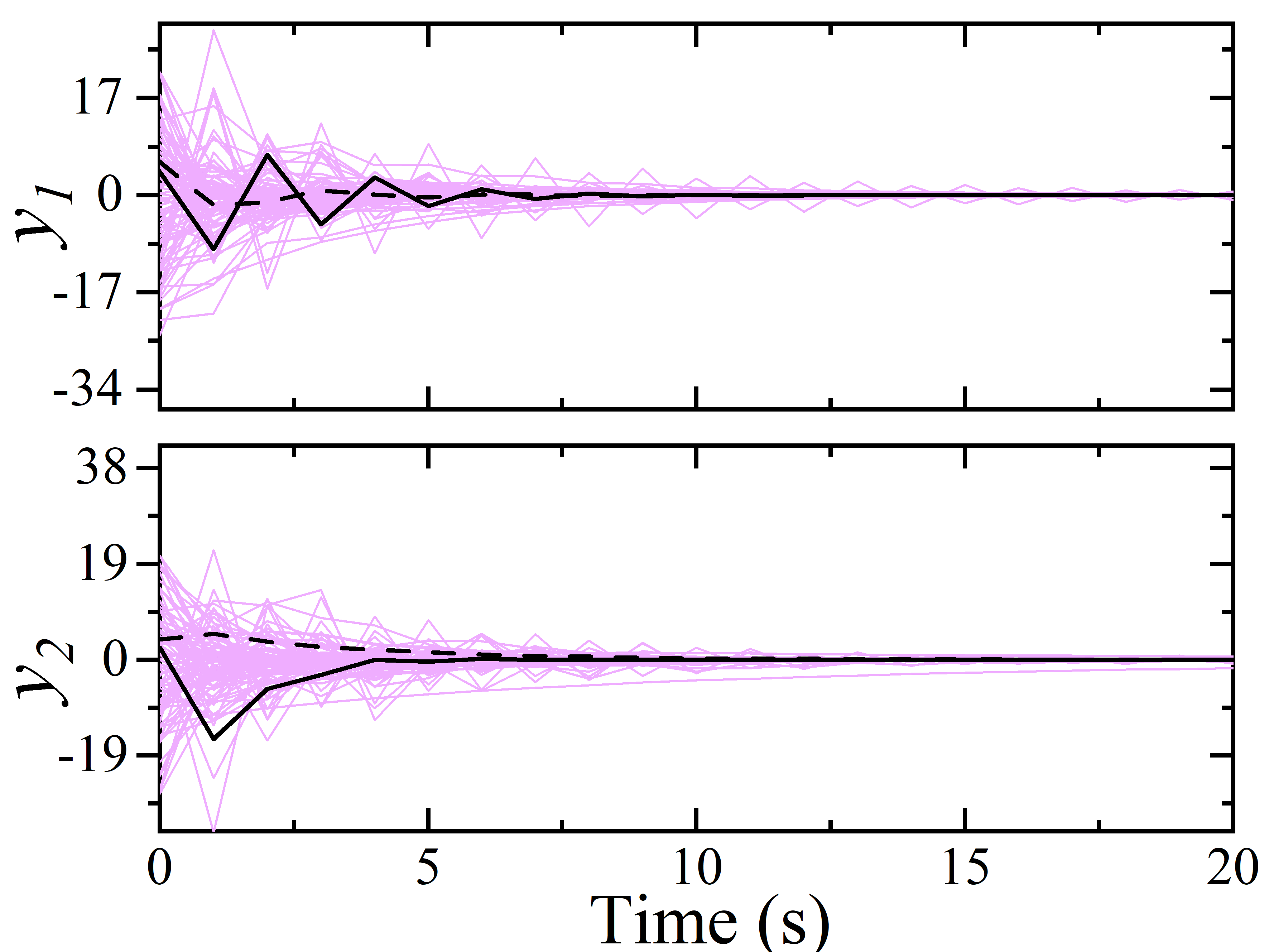}
            \caption{}
            \label{fig:2SubsOP}
        \end{subfigure}
        \begin{subfigure}[b]{0.23\textwidth}
            \includegraphics[width=\linewidth]{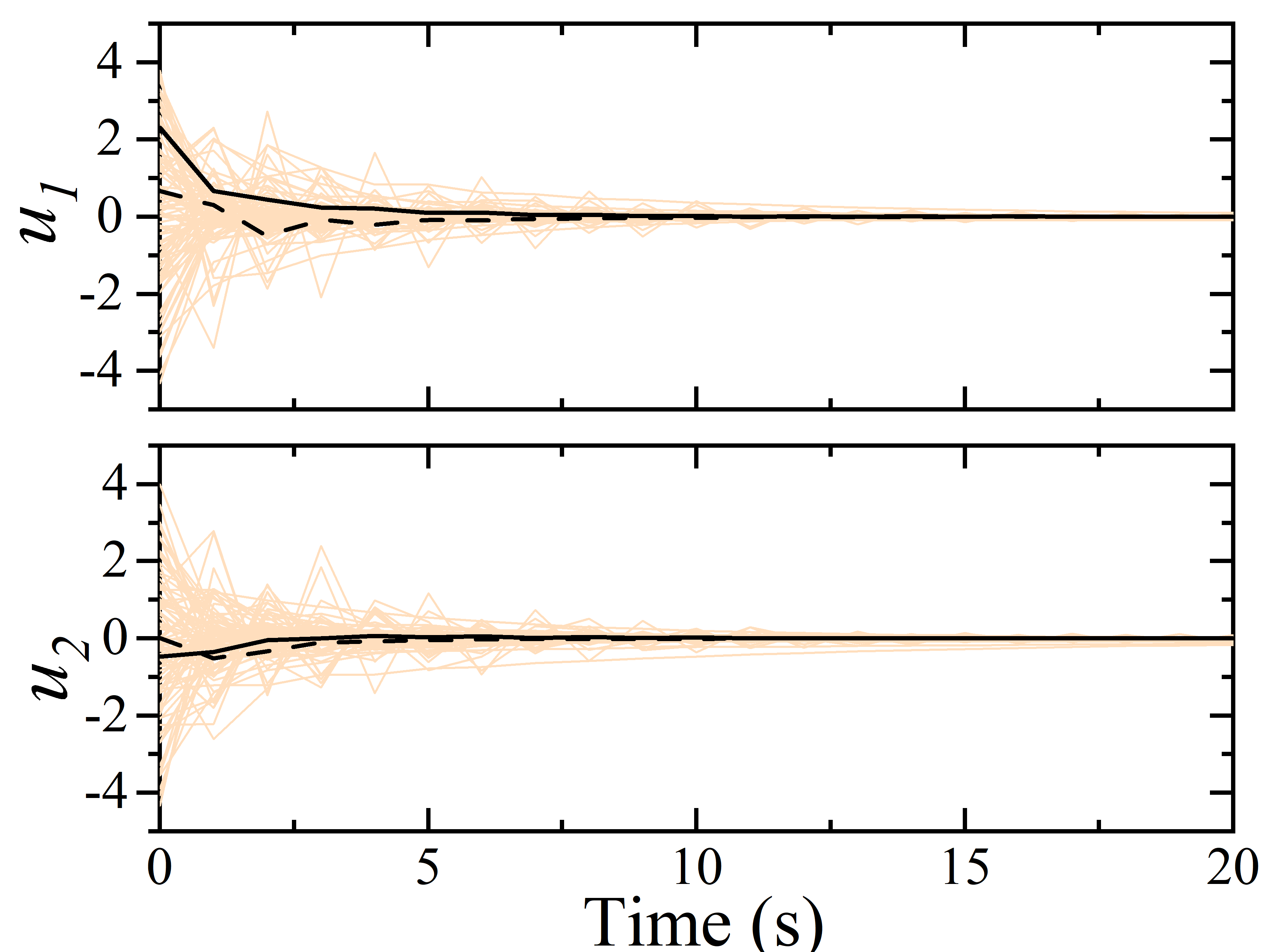}
            \caption{}
            \label{fig:2SubsIP}
        \end{subfigure}\\
  \caption{Controller performance for 2 subsystems. (a) Subsystem outputs. (b) Controller outputs. The black lines represent two random runs for visualization.}
  \label{fig:2SubsPerformances}
\end{figure}

The random plants were simulated for $100$ seconds with consistent random initial conditions for comparison across algorithms. To simplify results, the sum of all states for each subsystem is used as the output. For instance, in the case of $M = 2$ subsystems (Fig.~\ref{fig:2SubsPerformances}), each subsystem has 2 states $(\textbf{x}_1 \in \mathbb{R}^2, \textbf{x}_2 \in \mathbb{R}^2)$, and the outputs $y_1$ and $y_2$ are the element-wise sums of the state vectors. While this simplification may not be meaningful for real systems, it effectively demonstrates dynamic behavior and convergence under different control strategies. To ensure comparable control performance between DiMPC, I-mpDiMPC, IF-mpDiMPC, and FACET-DiMPC, a small error tolerance $\epsilon = 10^{-8}$ \cite{wegstein1958accelerating} was set with a maximum of $p_{max} = 100$ iterations. Results from all 100 runs for the case of 2 subsystems are shown in Fig.~\ref{fig:2SubsPerformances}. Due to the small tolerance, the outputs for all controllers are practically identical, achieving a centralized-like performance. 

For coupled systems, where the control inputs of one subsystem affect others, DiMPC and I-mpDiMPC require more intermediate iterations to converge, increasing communication and computational costs. In contrast, the proposed FACET-DiMPC approach can achieve centralized-like performance without iterations, directly reducing both communication overhead and computational load by solving linear equations instead of optimization problems.

\subsection{Communication Burden}
The maximum number of intermediate iterations per sample time step was recorded for each experiment, and the maximum and average values for DiMPC are shown in Table~\ref{tab: maxAndAvgIter}. As the number of subsystems increases, both the maximum and average iteration counts rise, consistent with findings by \cite{wang2022improved}. This highlights the increased computational effort required to reach a plantwide optimum using iterative methods.

\begin{table}
    \centering
     \caption{Maximum and average number of iterations under DiMPC.}
    \begin{tabular}{ccc}
    \hline
    \textbf{No. of Subsystems} & \textbf{Max Iterations} & \textbf{Aver. Iterations} \\
    \hline
     $2$ & $38$ & $19.99$  \\
    $3$ &  $62$ & $33.97$ \\
    $4$  &  $100$  &  $49.86$ \\
    $5$  & $100$ & $66.05$ \\
    \hline
    \end{tabular}
    \label{tab: maxAndAvgIter}
\end{table}

\begin{figure}[htbp]
    \centering
    \includegraphics[width=\columnwidth]{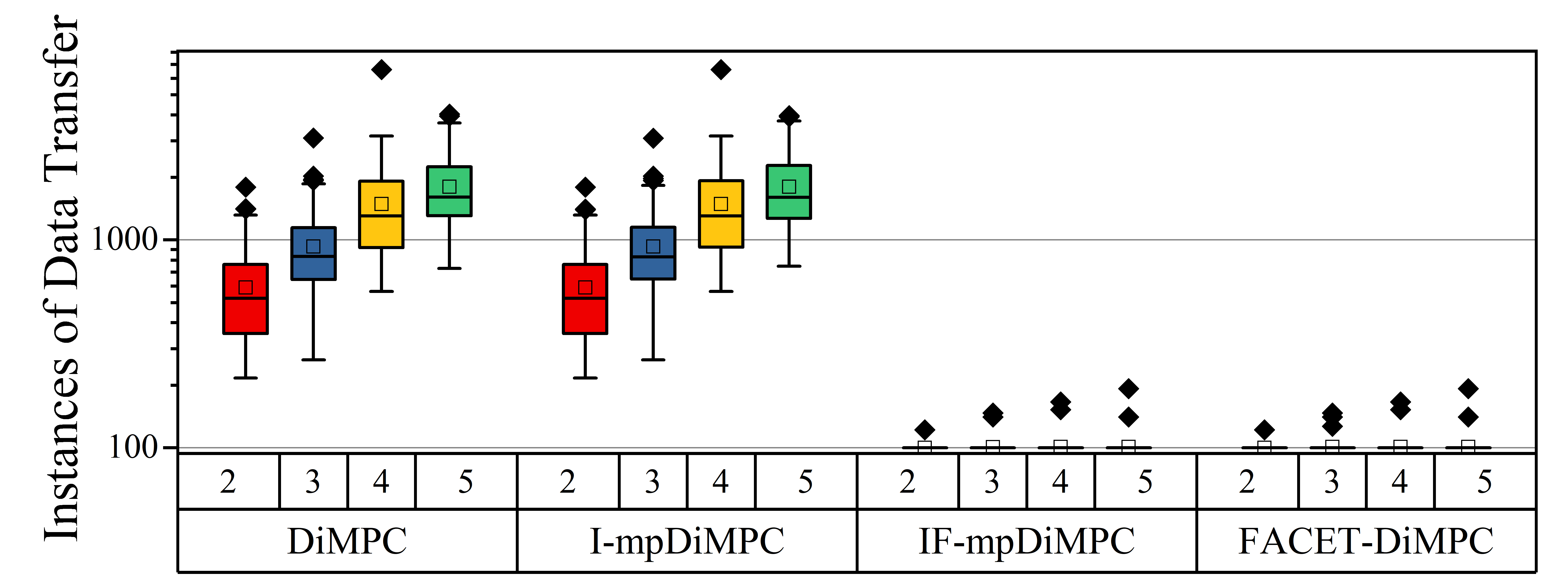}
    \caption{Number of data transfer instances for different DiMPC controllers during the entire simulation.}
    \label{fig:Comutimelog}
\end{figure}

To compare the communication burden, we analyzed the instances of data transfer between subsystems for all methods. Fig.~\ref{fig:Comutimelog} shows this comparison. The iteration-free methods (IF-mpDiMPC and FACET-DiMPC) significantly reduce data transfer instances, particularly as subsystems increase. Unlike iterative methods, which require multiple exchanges to converge, iteration-free methods only need near to single data transfer per time step (100 transfers for the entire 100-second simulation), regardless of subsystem count. This reduction in communication load decreases network traffic and improves system robustness, addressing one of the key motivations of this work.
\begin{figure}[htbp]

    \centering
    \includegraphics[width=0.95\columnwidth]{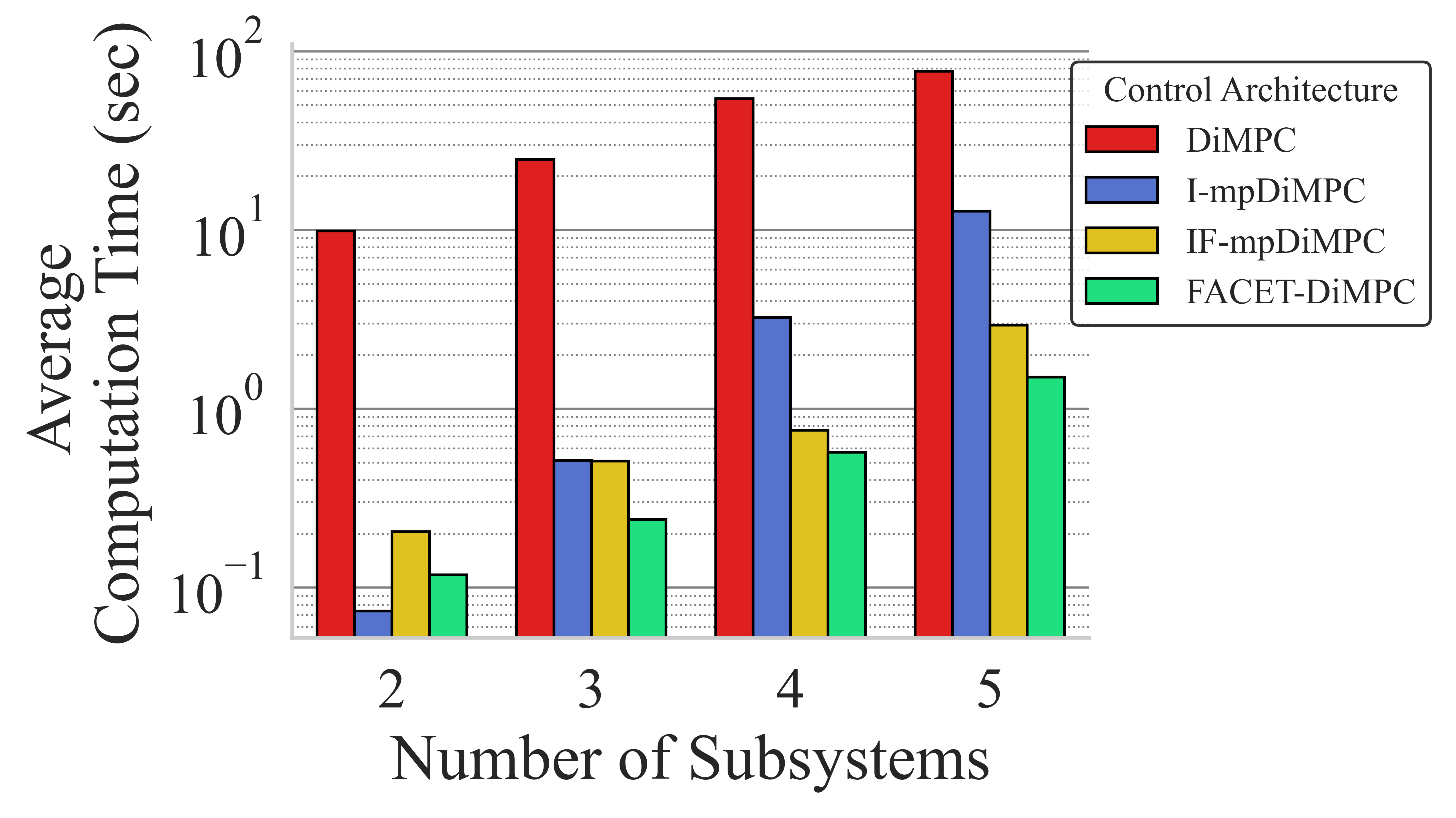}
    \caption{Average computation time (sec) of DiMPC control architectures over 100 plants. (Log-Scale)}
    \label{fig:MeanCompTims}
\end{figure}

\begin{figure}[htbp]
    \centering
    \includegraphics[width=\columnwidth]{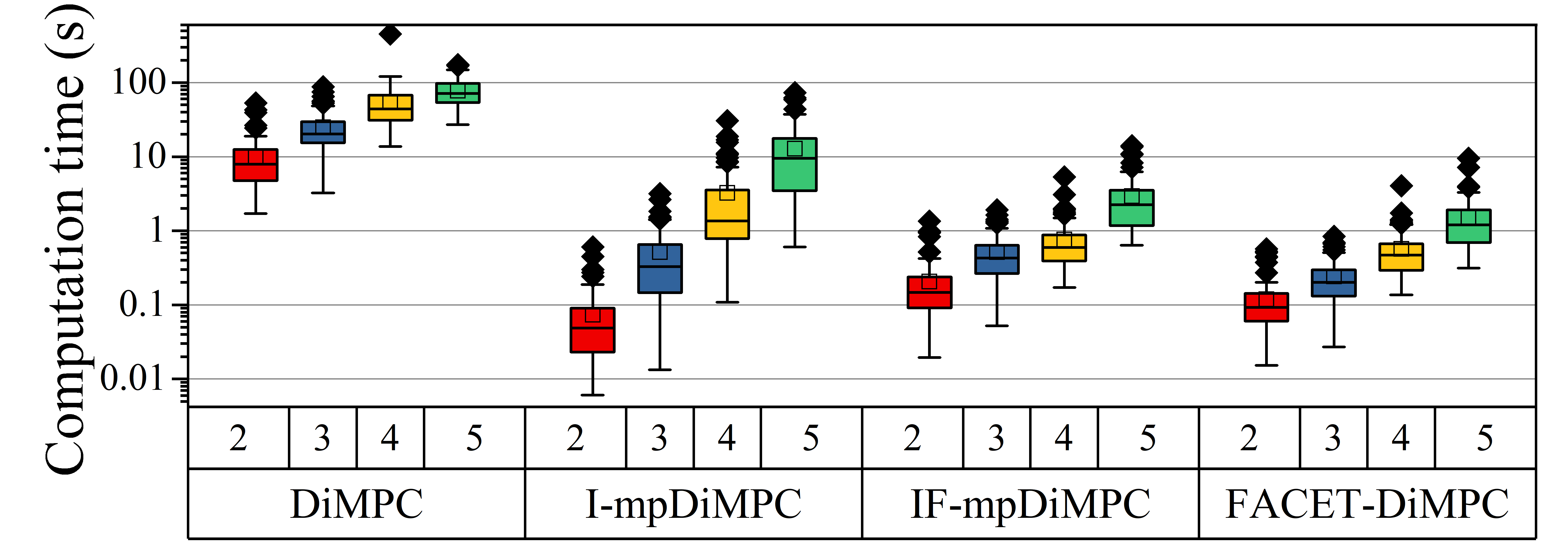}
    \caption{Computation times for different DiMPC controllers.}
    \label{fig:Comptimelog}
\end{figure}

\begin{figure}[htbp]
    \centering
    \includegraphics[width=\columnwidth]{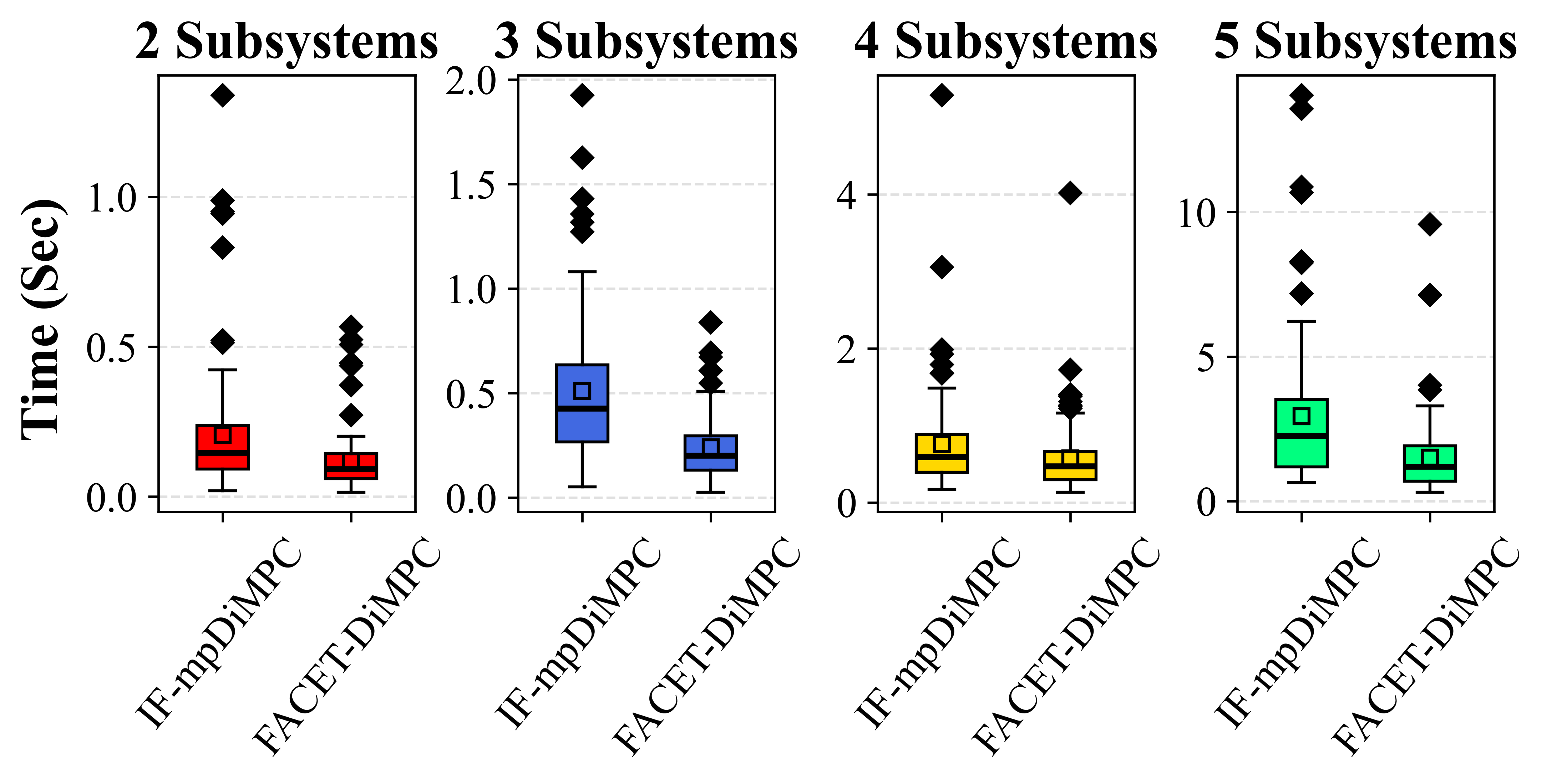}
    \caption{Computation times of IF-mpDiMPC and FACET-DiMPC controllers for different subsystems. (Linear-Scale)}
    \label{fig:Comptimelinear}
\end{figure}  
\subsection{Computation Time}

Figure~\ref{fig:Comptimelog} shows the computational performance of the controllers for a range of subsystem numbers. The computation time reflects the total time taken for all plant states to reach 0.02\% of their equilibrium value. Iteration-free DiMPC solves simultaneous equations at each time step, making its computation time scale linearly with simulation duration (100 s in this study). However, as the number of subsystems increases, the computation time for DiMPC grows exponentially, making it intractable for large systems and underscoring the need for efficient algorithms. As illustrated in Fig.~\ref{fig:Comptimelinear}, the proposed FACET-DiMPC method demonstrates a significant reduction in computation time compared to IF-mpDiMPC, an advantage that becomes increasingly pronounced as the number of subsystems grows from two to five. Figure~\ref{fig:MeanCompTims} compares the different approaches, showing that all mp-based methods have a significantly lower average computational time than classical DiMPC. While the iteration-free DiMPC reduces communication overhead, its higher computation times make it less practical for large-scale implementation. FACET-DiMPC, however, performs the best overall by effectively balancing computational efficiency and communication reduction. This superior performance is due to the restricted search space used in FACET-DiMPC, which limits the number of critical regions that need to be explored at each time step.

The results show that the FACET-DiMPC provides the best trade-off between computation time and communication load, while achieving similar control performance to the centralized methods. This makes FACET-DiMPC suitable for large-scale systems where communication delays and computational costs are critical constraints.

\section{Conclusions}

A novel iteration-free DiMPC method, FACET-DiMPC, was developed using mp programming while exploiting facet properties. The new method is compared against iterative and iteration-free counterparts. While the offline computational cost scales exponentially with the number of subsystems and states due to the increasing number of critical regions, a known limitation of mp-based approaches, FACET-DiMPC consistently achieved the fastest online solution times with minimal increase in communication overhead.  While current simulations were executed on a CPU, further enhancements to real-time scalability for large-scale systems could be achieved by leveraging GPU-based parallelization for concurrent region combinations evaluation, or by employing advanced point location techniques like binary search trees\cite{tondel2003evaluation} and machine learning based approximation\cite{karg2020efficient,airan2013linear}. Future work will apply the developed method on a relevant industrial case study, and extend its applicability to unstable and nonlinear systems.

\section*{Acknowledgment}

We thank Dr. Radhe Saini for the development of the IF-mpDiMPC code, that served as the foundation for this work.


\end{document}